\documentclass{article}
\usepackage{amsmath,amsthm,amsfonts,latexsym,amscd,amssymb,fancyhdr}
\newtheorem{theorem}{Theorem}[section]
\newtheorem{lemma}[theorem]{Lemma}
\newtheorem{corollary}[theorem]{Corollary}

\theoremstyle{definition}
\newtheorem{definition}[theorem]{Definition}

\theoremstyle{remark}

\newcommand{\N}{\mathbb{N}}

\newcommand{\Z}{\mathbb{Z}}

\newcommand{\n}{\noindent}
\begin{document}

\begin{center}
  {\large On Odd Order Nilpotent Group With Class $2$}
 \end{center}
\begin{center}
V. K. Jain\\
School of Mathematics, HRI, Allahabad, India\\
email: vkj@mri.ernet.in
\end{center}

\begin{abstract}
Let $G$ be an odd order nilpotent group with class $2$ and $e$ denotes the exponent of its commutator subgroup. Let $e=p_1^{r_1}p_2^{r_2}\ldots p_s^{r_s}$, where $p_i$'s are odd primes and $r_i$'s are non-negative integers. Then there are at least $r_1+r_2+\ldots +r_s$ non-isomorphic nilpotent groups with class two and the order of each of the group is equal to the order of $G$.
\end{abstract}



\section{Introduction and Statement of Main Result}
A very basic problem in the finite group theory is to calculate the number of non-isomorphic groups of a given order. One of the reduced version of this problem is for given a group with some group theoretic property, is there another non-isomorphic group (of equal order) with same property? In this paper, we construct a string of non-isomorphic groups from a given odd order nilpotent group with class $2$ (by a nilpotent group with class $2$, we mean that the group is non-abelian and the commutator subgroup of the group is contained in the center of the group).

Let $(G,.)$ be a group. Fix $n \in \N$, the set of natural number. Define a binary operation $\circ^n$ on $G$ as follows; for $x,y \in G$, $x\circ^n y:=y^{-n}xy^{n+1}$. Then equation of the type $X \circ^n a = b$, where $X$ is unknown and $a,b \in G$ has unique solution. But if $(G,.)$ is nilpotent group with class $2$, then $\circ^n $ is associative, that is $(G,\circ^n )$ is a group. Moreover, it can be observed that if $G$ is nilpotent group with class $2$, then $(G, \circ^n )$ is a nilpotent group with class at most $2$ and $x \circ^n y =[x,y]^nxy $ (using commutator identities).

\begin{definition}\label{vi}
Given a nilpotent group $(G,.) $ with class at most $2$, we define another group $(S_n(G),\circ ^n)$ as follows: $S_n(G)=G$, for $x,y \in S_n(G),~x \circ^n y: =[x,y]^nxy$. Further, we define a sequence $\{S_n^i(G)\}_i$ by induction. Let ${S_n^0(G)}= G$. If ${S_n^i(G)}$ is defined, then ${S_n^{i+1}(G)} = S_n(S_n^i(G))$.
\end{definition}

\begin{definition}
Let $(G,.)$ be a group of order $r$. An ordered $r-$ tuple \\$(n_1,n_2, \ldots , n_r)$ in $\N^r$ such that $1=n_1 \leq n_2 \leq \ldots \leq n_r$ is called {\em{order structure}} of $G$ if elements of $G$ can be put in a sequence $\{x_1, x_2, \ldots , x_r \}$ of length $r$ such that the order of $x_i$ in $G$ is $n_i, ~1\leq i \leq r$.
\end{definition}

For example, the order structure of $\Z_2 \times  \Z_4$ is $\{1,2,2,2,4,4,4,4\}$.
By the classification of finite abelian groups, it is clear that two abelian groups are isomorphic if and only if they have same order structure.

It is clear from the Definition \ref{vi} that the order of an element $x \in G$ is same as the order of $x \in S_n^i(G)$ for any $i$ and $n$. Thus the order structure of $S_n^i(G)$ is same as that of $G$. Also if $G$ is abelian, then $(S_n(G), \circ^n)=(G,.)$.

Let $r$ be an odd integer greater than $1$.
Let $r=p_1^{r_1}p_2^{r_2}\ldots p_s^{r_s}$,
where $p_j$'s are odd primes such that $p_1<p_2< \ldots <p_s$ and $r_j>0$. 
Let $X^r$ denotes a complete set of non-isomorphic nilpotent groups with class at most $2$ and order $r$. We say two elements of $X^r$ are related if both have same order structure.
 This relation is an equivalence relation. 
Let $X_1^r, X_2^r, \ldots , X_k^r$ be distinct equivalence classes of the set $X^r$ under this relation.
Let $G$ be a nilpotent group with class at most $2$ and order $r$. 
Identify $G$ with its unique isomorphic copy in $X^r$. 
We may assume that $G \in X_1^r$.
 
Since $G$ is nilpotent, $G=P_1P_2 \ldots P_s$ where $P_j$ is Sylow $p_j$-subgroup of $G$, $1 \leq j \leq s$. Take $n_j=\frac{p_j-1}{2}$ for all $1 \leq j \leq s$. Suppose that the exponent of the commutator subgroup of $P_j$ is $p_j^{t_j}$. We define the {\em{string of groups}} associated with $G$ to be a sequence $F_i(G)$ of length $\sum_{j=1}^{j=s}t_j +1$ in $X^r_1$ as 
 $F_i(G)$ $\cong $

\noindent 
$\left( \!\begin{array}{lr}
S^i_{{n_1}}(P_1)P_2P_3\ldots P_s & 0 \leq i \leq t_1   \\
S^{t_1}_{n_1}(P_1) S^{i-t_1}_{n_2}(P_2 )P_3 \ldots  P_s & t_1 < i \leq t_1+t_2  \\
\vdots & \vdots \\ S^{t_1}_{n_1}(P_1)S^{t_2}_{n_2}(P_2)   \ldots S^{t_{s-1}}_{n_{s-1}} (P_{s-1})S^{i-\sum_{j=1}^{j=s-1}t_j}_{n_s}(P_s)  & \sum_{j=1}^{j=s-1}t_j < i \leq \sum_{j=1}^{j=s}t_j
\end{array} \right)$.

\noindent For simplicity, we call the string of groups associated with $G$, {\em{the string of}} $G$. 
We call an element $G \in X_i^r$ { \em maximal} if there does not exist $H ~(\neq G) \in X_i^r$ such that $F_j(H) \cong G$ and call $G \in X_i^r$ {\em minimal} element  if $G$ is abelian.

Note that $X_i^r$ can have more than one maximal elements. For example, take $r=p^4$, where $p$ is an odd prime, then the following is a complete set of non-isomorphic nilpotent groups with class $2$ of order $r$ (see \cite[p. 145]{burn}).

\n (i) $A=\langle x,y |x^{p^3}, y^p, y^{-1}xy=x^{1+p^2}  \rangle $,

\n (ii) $B=\langle x,y,z |x^{p^2}, y^p, z^p, z^{-1}yz=yx^p, y^{-1}xy=x, z^{-1}xz=x \rangle $,

\n (iii) $  C=\langle x,y |x^{p^2},y^{p^2}, y^{-1}xy=x^{1+p} \rangle $,

\n (iv) $D=\langle x, y, z |x^{p^2}, y^p, z^p, z^{-1}xz=x^{1+p}, x^{-1}yx=y,z^{-1}yz=y \rangle $,

\n (v) $E= \langle x, y, z |x^{p^2}, y^p, z^p, z^{-1}xz=xy, y^{-1}xy=x,z^{-1}yz=y \rangle $,
 
\n(vi) $F= \langle x,y,z,a,| x^p,y^p,z^p,a^p, a^{-1}za=zx, a^{-1}ya=y, a^{-1}xa=x, z^{-1}yz=y, z^{-1}xz=x, y^{-1}xy=x \rangle $.

\n Then we may take $X^r=\{ A,B,C,D,E,F, \Z_{p^3} \times \Z_p, \Z_{p^2} \times \Z_p \times \Z_p, \Z_{p^2} \times \Z_{p^2}, \Z_p \times
\Z_p \times \Z_p \times \Z_p , \Z_{p^4} \}$, where $\Z_n$ denotes the cyclic group of order $n$. Equivalence classes of $X^r$ with respect to order structure relation are $X_1^r=\{ A, \Z_{p^3} \times \Z_p \}, ~ X_2^r=\{B,D,E, \Z_{p^2} \times \Z_p \times \Z_p \},~ X_3^r=\{ C, \Z_{p^2} \times \Z_{p^2} \} ~ X_4^r=\{ F, \Z_p \times \Z_p \times \Z_p \times \Z_p \}$ and $ X_5^r=\{\Z_{p^4}\} $.
In $X_2^r$, the commutator subgroups of $B,D$ and $E$ have exponent $p$, so the length of the string of each of these groups will be $2$. Thus $B,D$ and $E$ are three maximal elements in $X_2^r$.
The following is the main result of the paper.
\begin{theorem}\label{0}
{(a)} No two terms of the string of a group $G$ are isomorphic.\\
{(b)} The last term of the string of a group $G$ is an abelian group.\\
{(c)} The last terms of string of any two members of $X_i^r$ are same. 
\end{theorem} 

It is an interesting problem to determine all maximal elements in each $X^r_i$ for all $i$ because by knowing the maximal elements in $X^r_i$ for each $i$, one knows all elements of $X^r$.

\section{Basic Lemmas}
The results of this section has been proved using the commutator identities of a nilpotent group with class $2$. These identities are stated in the Lemma \ref{a}. We use $[x,y]$ for $x^{-1}y^{-1}xy$, where $x,y \in G$ and $Z(G)$ for the centre of the group $G$.

\begin{lemma} \label{a}
 Let $G$ be a nilpotent group with class $2$. Let $x,y,z \in G$. Then

\noindent (i) $[xy,z]=[x,z][y,z]$,

\noindent (ii) $[x^k,y]=[x,y^k]=[x,y]^k$, where $k$ is an integer.
\end{lemma}
\begin{proof}
 The first identity is easy to verify. The second can be proved by induction on $k$ and using (i).
\end{proof}

\begin{lemma}\label{1}
Let $G$ be a nilpotent group with class $2$ and $n\in \N$. Then $S_n^i(G)=S_{s(i)}(G)$, where $s(i)=\frac{(2n+1)^i-1}{2}$.
\end{lemma}
\begin{proof}
The proof is by induction on $i$. For $i=1$, $s(i)=n$, so identity is true. Suppose that identity is true for $i=m$. We will show that it is true for $i=m+1$.
Now $S_n^{m+1}(G)=S_n(S_n^m(G))=S_n(S_{s(m)}(G))$. Take $x,y \in S_n(S_{s(m)}(G))$. Then $x\circ ^n y=(x^{-1}\circ^{s(m)}y^{-1}\circ^{s(m)}x\circ^{s(m)}y)^n\circ^{s(m)}x\circ^{s(m)}y
=[x,y]^{(2s(m)+1)n+s(m)}xy$.
Now $2s(m)n+n+s(m)=s(m)(2n+1)+n=\frac{(2n+1)^m-1}{2}(2n+1)+n=\frac{(2n+1)^{m+1}-1}{2}=s(m+1)$. This proves the lemma.
\end{proof}
\begin{lemma}\label{2}
Let $G$ be a nilpotent group with class two and $|G|=p^k$ where $p$ is an odd prime. Then for $n=\frac{p-1}{2}$, we have

\noindent (i) $x\in Z(S_n^i(G)) \iff x^{p^i} \in Z(G)$.

\noindent (ii) $x\in Z(S_n^{i+1}(G)) \iff x^p \in Z(S_n^i(G))$.

\noindent (iii) $Z(S_n^i(G))$ is a normal subgroup of $S_n^{i+1}(G)$.
\end{lemma}
\begin{proof}
$(i)$ Let $x\in Z(S_n^i(G))$ and $y\in S_n^i(G)$.
By Lemma \ref{1}, $S_n^{i}(G)=S_{s(i)}(G)$, so 

$ 
\begin{array}{lccl}
& x \circ^{s(i)} y&=& y\circ^{s(i)}x \\
\iff & x\circ^{s(i)}y\circ^{s(i)}x^{-1} \circ^{s(i)}y^{-1}&=&1\\
\iff & [x,y]^{2s(i)+1}&=&1 ~~~ (\text{Lemma}~ \ref{a})\\
\iff & [x,y]^{p^i}&=&1~~~(\text{for all}~ y\in S_n^i(G))\\
\iff & [x^{p^i},y]&=&1~~~(\text{Lemma}~ \ref{a}(ii))\\
\iff & x^{p^i} &\in & Z(G).\\
\end{array}
$

This proves the (i). Second part follows from (i).

\noindent $(iii)$ From (ii) it is clear that $Z(S_n^i(G))$ is a subset of $S_n^{i+1}(G)$. 
By Lemma \ref{1}, $S_n^i(G)=S_{s(i)}(G)$. 
Let $x,y \in Z(S^i_n(G))$. Then $x \circ ^{s(i)} y=y \circ ^{s(i)} x$. 
This implies $[x,y]^{s(i)}xy=[y,x]^{s(i)}yx$, that is $[x,y]^{2s(i)+1}=1$. Now $s(i+1)=\frac{(2n+1)^{i+1}-1}{2}=\frac{(2n+1)^i(2n+1)-1}{2}=(2s(i)+1)n+s(i)$. 
Now since $x,y \in S^{i+1}_n(G)$ and $S^{i+1}_n(G)=S_{s(i+1)}(G)$, so 

$\begin{array}{lcl}
x \circ ^{s(i+1)}y &=& [x,y]^{s(i+1)}xy\\
&=& [x,y]^{(2s(i)+1)n+s(i)}xy \\
&=& [x,y]^{s(i)}xy ~(\text{for} ~~[x,y]^{2s(i)+1}=1)\\
&=& x \circ^{s(i)} y.
\end{array}$

\noindent   This implies $Z(S_n^i(G))$ is a subgroup of $S^{i+1}_n(G)$. By (ii), $Z(S_n^i(G))\subseteq Z(S_n^{i+1}(G))$. Thus $Z(S_n^i(G))$ is a normal subgroup of $S_n^{i+1}(G)$.
\end{proof}

\section{Proof of the Theorem}
\begin{proof}[Proof of the Theorem \ref{0}]
 It is sufficient to prove the result for $p$-group, where $p$ is an odd prime. Let $p$ be an odd prime and $G$ be a nilpotent group of order $r=p^k$ with class $2$. Let the exponent of the commutator subgroup of $G$ be $p^t$. Assume that $G \in X_1^r$ (see Section 1). Identify $F_i(G)$ with $S^i_{\frac{p-1}{2}}(G), 0 \leq i \leq t$. Take $n=\frac{p-1}{2}$. \\
\textbf{(a)}
 To prove it we will show that centers of two distinct terms of the string of $G$ are not equal. Fix $i \in \{0, \ldots , t-1 \}$. By Lemma \ref{2} (iii), $Z(S_n^i(G))$ is a normal subgroup of $S_n^{i+1}(G)$. Let $H=S^{i+1}_n(G)/Z(S^i_{n}(G))$. 
Take $\Omega_1(H)= \langle \{ \bar{x} \in H| \bar{x}^p=\bar{1} \in H \} \rangle$, where $\bar{1}$ denotes the identity element of $H$. 
But since $H$ is nilpotent, $\Omega_1(H) \neq \{\bar{1}\}$.
This implies that there exists $x \in S^{i+1}_n(G)$ such that $x \not \in Z(S^i_{n}(G))$ and $x^p \in Z(S^i_{n}(G))$. Thus $Z(S^i_{n}(G)) \subsetneq Z(S^{i+1}_{n}(G))$ for all $0 \leq i \leq t-1$.
That is, $Z(G) \subsetneq Z(S_n(G)) \subsetneq Z(S^2_n(G)) \subsetneq \ldots \subsetneq Z(S^t_n(G))$. Hence distinct terms of the string of $G$ are non-isomorphic groups. 

\n \textbf{(b)}
The last term of the string of $G$ is $S_n^t(G)$.
By Lemma \ref{1}, $S_n^t(P)=S_{s(t)}(G)$ where $s(t)=\frac{(2n+1)^t-1}{2}$. But $n= \frac{p-1}{2}$, so $s(t)=\frac{p^t-1}{2}$.
Now let $x, y \in S_n^t(P)$. Then 

$\begin{array}{lcl}
  x \circ^{s(t)}y  &=& [x,y]^{s(t)}xy\\
&=& [x,y]^{\frac{p^t-1}{2}}xy\\
&=& [x,y]^{p^t-\frac{p^t+1}{2}}xy\\
&=& [x,y]^{- \frac{p^t+1}{2}}xy ~~~~(\text{for}~p^t ~\text{is the exponent of the commutator})\\
&=&[ y,x]^{\frac{p^t+1}{2}}xyx^{-1}y^{-1}yx~~~~~~(\text{for}~[x,y]^{-1}=[y,x])\\
&=& [y,x]^{\frac{p^t+1}{2}}[x^{-1},y^{-1}]yx\\
&=& [y,x]^{\frac{p^t+1}{2}} [y,x]^{-1}yx ~~~~(\text{Lemma \ref{a}}(ii)) \\
&=& [y,x]^{\frac{p^t-1}{2}}yx\\
&=& y \circ ^{s(t)} x.
 \end{array}$

So, $S_n^t(P)$ is an abelian group.

\n \textbf{(c)} Let $G, H \in X_i^r$. Since the last term of the string is an abelian group and for a given order structure there is unique abelian group, so last term of $G$ and $H$ are same.
\end{proof}
The following is an easy consequence of the above Theorem.

\begin{corollary}
Let $(G,.)$ be an odd order nilpotent group with class two such that the exponent of the commutator subgroup is $p_1^{r_1}p_2^{r_2}\ldots p_s^{r_s}$. Then there are at least $r_1+r_2+\ldots +r_s$ non-isomorphic nilpotent groups with class two and the order structure of each of these groups is same as that of $G$. 
\end{corollary}
\begin{proof}
The length of the string of $G$ is $1+ r_1+r_2+\ldots +r_s$ and each term of string has same order structure. By Theorem \ref{0}, the terms of the string are non-isomorphic groups and last term of the string is abelian. Thus the string of $G$ except last term gives us the required non-isomorphic groups.
\end{proof}

\end{document}